\documentclass{amsart}
\usepackage{amsmath}
\usepackage{amssymb}
\usepackage{graphicx}
\usepackage{tikz}
\usepackage{enumitem}
\usepackage{comment}
\usepackage{bm}
\usepackage{soul}
\usepackage{fullpage}
\usepackage{amsaddr}

\newcommand{\Tb}{\mathbf{T}}
\newcommand{\etab}{\boldsymbol{\eta}}
\newcommand{\xib}{\boldsymbol{\xi}}
\newcommand{\alphab}{{\boldsymbol{\alpha}}}
\newcommand{\nb}{\mathbf{n}}
\newcommand{\lno}{\left \|}
\newcommand{\rno}{\right \|}
\newcommand{\res}[1]{\left.{#1} \right|}
\newcommand{\ex}{{\etab(\xib)}}
\newcommand{\EE}{\mathcal E}
\newcommand{\lp}{\left(}
\newcommand{\rp}{\right)}

\def\xib{{\bm\xi}}

\def\etab{{\bm\eta}}
\def\alphab{{\bm\alpha}}

\newtheorem{lemma}{Lemma}
\newtheorem{theorem}{Theorem}

\title{Enforcing Neumann Boundary Conditions with Polynomial Extension Operators to Acheive Optimal Convergence Rates on Polytopial Meshes in the Finite Element Method}
\author{James Cheung}
\address{Millennium Space Systems, A Boeing Company. 2265 E. El Segundo Blvd, El Segundo, CA.}

\begin{document}
\maketitle

\begin{abstract}
In \cite{cheung2019optimally}, the authors presented two finite element methods for approximating second order boundary value problems on polytopial meshes with optimal accuracy without having to utilize curvilinear mappings. This was done by enforcing the boundary conditions through judiciously chosen polynomial extension operators. The $H^1$ error estimates were proven to be optimal for the solutions of both the Dirichlet and Neumann boundary value problems. It was also proven that the Dirichlet problem approximation converges optimally in $L^2$. However, optimality of the Neumann approximation in the $L^2$ norm was left as an open problem. In this work, we seek to close this problem by presenting new analysis that proves optimal error estimates for the Neumann approximation in the $W^1_\infty$ and $L^2$ norms.
\end{abstract}
\section{Introduction}
The purpose of this note is to derive optimal error estimates for the polynomial extension finite element method (PE-FEM) described in \cite{cheung2019optimally} for approximating elliptic boundary value problems with Neumann conditions. This manuscript is very much an appendix to our previous work. As such, we highly suggest that the reader refer to that work, especially since we will not redefine our notation here for the sake of brevity. Our analysis will be structured in the following manner: We present new technical lemmas involving the Averaged Taylor series in \S2, then in \S3 we move on to prove well-posedness and derive optimal error estimates for the numerical solution in the $W^1_\infty(\Omega)$ norm, using this new result we are finally able to derive optimal error estimates for the numerical solution in the $L^2(\Omega)$ norm in \S4. We then discuss our results and propose additional future research in \S5.

\section{Additional Results for Averaged Taylor Polynomials}
We begin our analysis by deriving some technical results for the Averaged Taylor polynomials in the $W^m_\infty(\Omega_h)$ setting. 
\begin{lemma} \label{lemma: Taylor lemma 1}
Let $v\in L^\infty(\Omega_h)$, then
$$
	\lno\left.T^k_h(v) \right|_{\etab(\xib)}\rno_{L^\infty(\Omega_h)} \leq C\|v\|_{L^\infty(\Omega_h)}.
$$
\end{lemma}
\begin{proof}
	Using a scaling argument on \cite[Corollary 4.1.15]{brenner2008mathematical}, we immediately have that
	$$
		\begin{aligned}
			\lno \res{T^k_h(v)}_{\etab(\xib)} \rno_{L^\infty(\EE_i \cap S^{i,\ell})} &\leq C\delta_h^{-d} \|v\|_{L^1(\sigma^{i,\ell})}\\
			&\leq C\|v\|_{L^\infty(\sigma^{i,\ell})}
		\end{aligned}
	$$
	We conclude by seeing that
	$$
		\lno \res{T^k_h(v)}_\ex \rno_{L^{\infty}(\Gamma_h)} = \max_{i,\ell} \lno \res{T^k_h (v)}_\ex \rno_{L^\infty(\EE_i \cap S^{i,\ell})}.
	$$
\end{proof}

\begin{lemma} \label{lemma: Taylor lemma 2}
Let $v \in W^{k+1}_\infty(\Omega_h)$. Then for an integer $m<k+1$, the following is satisfied
$$
	\lno \res{T^k_h(v)}_\ex - v \rno_{W^m_\infty(\Gamma_h)} \leq C \delta_h^{k+1-m} |v|_{W^{k+1}_\infty(\Omega_h)}.
$$
\end{lemma}
\begin{proof}
	Using \cite[Proposition 4.3.2]{brenner2008mathematical}, we have that
	$$
	\begin{aligned}
		\lno \res{T^k_h(v)}_\ex - v \rno_{W^m_\infty(\EE_i\cap S^{i,\ell})} &= \lno T^k_h(v) - v \rno_{W^m_\infty(\etab(\EE_i)\cap S^{i,\ell})} \\ 
		&\leq \lno T^k_h(v) - v \rno_{W^m_\infty(S^{i,\ell})} \\ 
		&\leq C \delta_h^{k+1-m} |v|_{W^{k+1}_\infty(S^{i,\ell})}.
	\end{aligned}
	$$
	Taking the maximum over all $i,\ell \in \mathbb N$ concludes the proof. 
\end{proof}

\begin{lemma}\label{lemma: Taylor lemma 3}
Let $v\in\overline V^k_h$, then
$$
	\lno \res{T^{k',k}_h(v)}_\ex \rno_{L^\infty(\Gamma_h)} \leq C \sum_{\alphab |=1}^k h^{-|\alphab|}\delta_h^{|\alphab|} \|v\|_{L^\infty(\Omega_h)}
$$
\end{lemma}
\begin{proof}
	From the definition of $T^{k',k}_h(v)$, we see directly that
	$$
		\begin{aligned}
		\lno \res{T^{k',k}_h(v)}_\ex \rno_{L^\infty(\Gamma_h)} 
		&\leq \sum_{|\alphab|=1}^k \frac{\delta_h^{|\alphab|}}{|\alphab| !} \lno D^\alphab v \rno_{L^\infty(\Gamma_h)}\\
		&\leq \sum_{|\alphab|=1}^k \frac{\delta_h^{|\alphab|}}{|\alphab| !} \lno D^\alphab v \rno_{L^\infty(\Omega_h)}\\
		&\leq \sum_{|\alphab|=1}^k \frac{\delta_h^{|\alphab|} h^{-|\alphab|} }{|\alphab| !} \lno v \rno_{L^\infty(\Omega_h)},
		\end{aligned}
	$$
	after applying the inverse inequality.
\end{proof}

With these technical lemmas derived, we are now ready to prove that the solution of the Neumann approximation presented in \cite{cheung2019optimally} is well-posed and optimal in $W^1_\infty(\Omega)$.

\section{Well-Posedness and Error Estimates in $W^1_\infty(\Omega_h)$}
In this section, we determine that the solution $u_h \in V^k_h$ is bounded in $W^1_\infty(\Omega_h)$. Additionally, we demonstrate that the error estimate is optimal in the same norm. The analysis presented here is remarkably standard since the perturbations incurred by the extensions in the discrete bilinear form are not large enough to cause stability loss. Additionally, Strang's Lemma arguments are also used to handle the error incurred by the nonconforming perturbations in the bilinear form. 
\begin{theorem} \label{theorem: well-posedness theorem}
	Assume that $\delta_h\sim\mathcal O(h^2)$ and that $h \in \mathbb R_+$ is sufficiently small. Further assume that $\widetilde p,\widetilde q \in C^{k+1}(\Omega_h)$. Then there exist positive constants $\gamma, M \in \mathbb R_+$ such that
	\begin{equation} \label{eqn: well-posedness 1}
		\sup_{v_h\in V^k_h \cap W^1_1(\Omega_h)} \frac{B_{h,N}(w_h, v_h)}{ \|v_h\|_{W^1_1(\Omega_h)}} \geq \gamma\|w_h\|_{W^1_\infty(\Omega_h)}
	\end{equation}
	and
	\begin{equation} \label{eqn: well-posedness 2}
		B_{h,N}(w_h, v_h) \leq M\|w_h\|_{W^1_\infty(\Omega_h)} \|v_h\|_{W^1_1(\Omega_h)}
	\end{equation}
	for all $u_h, v_h \in V^k_h$.
\end{theorem}
\begin{proof}
	We begin by recalling that 
	$$
		B_{h,N}(w_h,v_h) = \int_{\Omega_h} \lp \widetilde p\nabla w_h \cdot \nabla v_h + \tilde q w_h v_h\rp dx + \left<\widetilde p\circ \ex \res{\Tb^{k-1}_h(\nabla w_h)}_\ex\cdot \nb - \widetilde p(\xib) \nabla w_h\cdot\nb_h, v_h \right>_{\Gamma_h}.
	$$
	We then see that
	$$
		\begin{aligned}
			B_{h,N}(w_h,v_h) \geq \int_{\Omega_h} \lp \widetilde p\nabla w_h \cdot \nabla v_h + \widetilde q w_hv_h \rp dx - Ch|u|_{W^1_\infty(\Omega_h)}\|v_h\|_{W^1_1(\Omega_h)},
		\end{aligned}
	$$
	after applying Lemma \ref{lemma: Taylor lemma 3} and the trace inequality \cite[Theorem 1.6.6]{brenner2008mathematical}. Dividing both sides with $\| v_h \|_{W^1_1(\Omega_h)}$ allows us to see that
	$$
		\sup_{v_h\in V^k_h \cap W^1_1(\Omega_h)} \frac{B_{h,N}(w_h, v_h)}{ \|v_h\|_{W^1_1(\Omega_h)}} \geq 
		\sup_{v_h\in V^k_h \cap W^1_1(\Omega_h)} \frac{\int_{\Omega_h} \lp \widetilde p\nabla w_h \cdot \nabla v_h + \widetilde q w_hv_h \rp dx}{\|v_h\|_{W^1_1(\Omega_h)}} - Ch|u|_{W^1_\infty(\Omega_h)}.
	$$
	Now, choosing $v_h = w_h$ allows us to see that $\sup_{v_h\in V^k_h \cap W^1_1(\Omega_h)}  \frac{\int_{\Omega_h} \lp \widetilde p\nabla w_h \cdot \nabla v_h + \widetilde q w_hv_h \rp dx}{\|v_h\|_{W^1_1(\Omega_h)}} > 0$. As such, we can choose a constant $C_{p,q} > 0$ such that
	$$
	\begin{aligned}
	\sup_{v_h\in V^k_h \cap W^1_1(\Omega_h)}  \frac{\int_{\Omega_h} \lp \widetilde p\nabla w_h \cdot \nabla v_h + \widetilde q w_hv_h \rp dx}{\|v_h\|_{W^1_1(\Omega_h)}} &\geq C_{p,q} \sup_{v_h\in V^k_h \cap W^1_1(\Omega_h)}  \frac{\int_{\Omega_h} \lp \nabla w_h \cdot \nabla v_h +w_hv_h \rp dx}{\|v_h\|_{W^1_1(\Omega_h)}} \\
	&= C_{p,q} \|w_h\|_{W^1_\infty(\Omega_h)}.
	\end{aligned}
	$$
	Therefore,
	$$
		\sup_{v_h\in V^k_h \cap W^1_1(\Omega_h)} \frac{B_{h,N}(w_h, v_h)}{ \|v_h\|_{W^1_1(\Omega_h)}} \geq C_{p,q} \|w_h\|_{W^1_\infty{\Omega_h}} - Ch \|u_h \|_{W^1_\infty(\Omega_h)}.
	$$
	And thus, we have that \eqref{eqn: well-posedness 1} is satisfied. 

	Using Lemma \ref{lemma: Taylor lemma 1}, H\"older's inequality, and the trace inequality allows us to derive \eqref{eqn: well-posedness 2}. 
\end{proof}

\begin{theorem} \label{theorem: max derivative optimality}
	Let $u_h \in V^k_h$ satisfy \cite[Equation 26]{cheung2019optimally}. Assume that $u \in W^{k+1}_{\infty}(\Omega)$ and that $f\in W^{k-1}_{\infty}(\Omega)$. Furthermore, let $\widetilde u \in W^{k+1}_\infty(\Omega_h)$ and $\widehat f, \widetilde f \in W^{k-1}_{\infty}(\Omega_h)$ be extensions of $u$ and $f$ respectively from $\Omega$ to $\Omega_h$. We then have that
	$$
		\|\widetilde u - u_h \|_{W^1_\infty(\Omega_h)} \leq Ch^{k} \lp |u|_{W^{k+1}_\infty(\Omega)} + \|f\|_{W^{k-1}_\infty(\Omega)} \rp
	$$
	under the conditions specified in Theorem \ref{theorem: well-posedness theorem}.
\end{theorem}
\begin{proof}
	Let $u_I \in V^k_h$ be the piecewise polynomial interpolant of $\widetilde u \in W^{k+1}_{\infty}(\Omega_h)$ defined on $\Omega_h$. From \cite[Theorem 4.2.20]{brenner2008mathematical} and the Stein extension theorem \cite{adams2003sobolev}, we have that
	\begin{equation} \label{eqn: interpolation error}
		\lno \widetilde u - u_I \rno_{W^1_\infty(\Omega_h)} \leq Ch^{k}|u|_{W^{k+1}_\infty(\Omega)}.
	\end{equation}

	We begin the analysis by seeing that
	$$
	\begin{aligned}
		\|u_I - u_h \|_{W^1_\infty(\Omega_h)} &\leq \sup_{v_h \in V^k_h\cap W^1_1(\Omega_h)} \frac{B_{h,N}(u_I - u_h, v_h)}{\|v_h\|_{W^1_1(\Omega_h)}} \\
		&= \sup_{v_h \in V^k_h\cap W^1_1(\Omega_h)} \frac{B_{h,N}(u_I - \widetilde u, v_h) + B_{h,N}(\widetilde u - u_h, v_h)}{\|v_h\|_{W^1_1(\Omega_h)}} \\
		&\leq M \|\widetilde u - u_h\|_{W^1_{\infty}(\Omega_h)} + \sup_{v_h \in V^k_h\cap W^1_1(\Omega_h)} \frac{\left<\widehat f - \widetilde f, v_h \right>_{\Omega_h} - \left< \widetilde p \circ \ex \mathbf{R}^{k-1}\res{(\nabla \widetilde u)}_\ex \cdot \nb, v_h \right>_{\Gamma_h}}{\|v_h\|_{W^1_1(\Omega_h)}}  \\
		&\leq Ch^k |u|_{W^{k+1}_\infty(\Omega)} + C\delta_h^{k-1}|f|_{W^{k-1}_\infty(\Omega)} + C\delta_h^{k-1}|u|_{W^{k+1}_\infty(\Omega)} \\
		&\leq Ch^k |u|_{W^{k+1}_\infty(\Omega)}
	\end{aligned}
	$$
	after utilizing \eqref{eqn: interpolation error}, Lemma \ref{lemma: Taylor lemma 2}, and seeing that
	$$
		B_{h,N}(\widetilde u, v_h) = \left< \widehat f, v_h \right>_{\Omega_h} + \left< g_N\circ\ex - \widetilde p\circ\ex \mathbf R^{k-1}\res{(\nabla u)}_\ex \cdot \nb, v_h \right>_{\Gamma_h}.
	$$

	The proof is completed by seeing that
	$$
		\|\widetilde u - u_h \|_{W^1_\infty(\Omega_h)} \leq \|u - u_I\|_{W^1_\infty(\Omega_h)} + \|u_I - u_h\|_{W^1_\infty(\Omega_h)},
	$$
	and applying the above bound along with \eqref{eqn: interpolation error}.
\end{proof}

\section{Error Estimates in $L^2(\Omega_h)$}
We are now ready to derive the optimal error estimates for the solution of the Neumann approximation in the $L^2(\Omega)$ norm. The analysis begins by estimating the nonconformity error. This nonconformity error will then be used in the following duality argument to bound the terms in the discrete problem that are not orthogonal in the Galerkin sense with respect to the continuous bilinear form. 

\subsection{Nonconformity Error}
Let us begin the derivation of the $L^2(\Omega_h)$ error bound by analyzing the nonconformity induced by $B_{h,N}(\cdot,\cdot):= H^1(\Omega_h)\times V^k_h \rightarrow \mathbb R_+$. 
\begin{lemma} \label{lemma: nonconformity error lemma}
Assume that all the conditions in Theorem \ref{theorem: max derivative optimality} hold, then the following is satisfied
$$
	N_h(\widetilde u - u_h, v_h) \leq Ch^{k+1}\lp |u|_{W^{k+1}_\infty(\Omega)} + |f|_{W^{k-1}_\infty(\Omega)} \rp \|v_h\|_{1,\Omega_h}
$$
for all $V^k_h \cap W^1_1(\Omega_h)$.
\end{lemma}
\begin{proof}
	Notice that
	$$
		B_{h,N}(\widetilde u, v_h) = \left< \widetilde f, v_h \right>_{\Omega_h} + \left< g_N\circ\ex - \widetilde p\circ\ex \res{\mathbf R^{k-1}_h(\nabla \widetilde u)}_\ex \cdot \nb, v_h \right>_{\Gamma_h} \quad \forall v \in H^1(\Omega_h).
	$$
	Taking the difference with \cite[Equation (24)]{cheung2019optimally} yields
	$$
		B_{h,N}(\widetilde u - u_h, v_h) = \left<\widehat f - \widetilde f, v_h \right>_{\Omega_h} - \left<\widetilde p\circ\ex \res{\mathbf R^{k-1}_h(\nabla \widetilde u)}_\ex \cdot \nb, v_h \right>_{\Gamma_h} \quad \forall v_h \in V^k_h.
	$$
	Let us define $e_h:=\widetilde u - u_h$, then from the definition of $B_{h,N}(\cdot,\cdot)$ (See \cite[Equation (25)]{cheung2019optimally}), we have that
	\begin{equation} \label{eqn: nonconformity inequality}
	\begin{aligned}
		N_h(e_h,v_h) &= \left< \widehat f - \widetilde f, v_h \right>_{\Omega_h} - \left<\widetilde p\circ\ex \res{\mathbf R_h^{k-1}(\nabla \widetilde u)}_\ex\cdot \nb,v_h \right>_{\Gamma_h} \\
					&\quad - \left< \widetilde p\circ\ex \res{\Tb^{k-1}_h(\nabla e_h)}_\ex\cdot \nb, v_h \right>_{\Gamma_h}
						   + \left< \widetilde p(\xib) \nabla e_h\cdot \nb_h,v_h\right>_{\Gamma_h} \\
				   &= \left< \widehat f - \widetilde f, v_h \right>_{\Omega_h} - \left<\widetilde p\circ\ex \res{\mathbf R_h^{k-1}(\nabla \widetilde u)}_\ex\cdot \nb,v_h \right>_{\Gamma_h} \\
				   &\quad + \left< \widetilde p\circ\ex \nabla e_h\cdot(\nb-\nb_h), v_h \right>_{\Gamma_h}
				   		  + \left< \widetilde p(\xib) \res{\Tb_h^{1,k-1}(\nabla e_h)}_\ex\cdot \nb, v_h\right>_{\Gamma_h} \\
				   &\quad + \left< (\widetilde p(\xib) - \widetilde p\circ\ex)\res{\Tb^{k-1}_h(\nabla e_h)}_{\ex}\cdot\nb,v_h\right>_{\Gamma_h}.
	\end{aligned}
	\end{equation}
	We will now analyze each of the terms on the right hand side of \eqref{eqn: nonconformity inequality} seperately. 

	Let $\Omega_{diff}^h:=\Omega_h \setminus (\Omega \cap \Omega_h)$, then we have by the definition of the extension that
	$$
		\begin{aligned}
			\left< \widetilde f - \widehat f, v_h \right>_{\Omega_h} 
			&= \left< \widetilde f - \widehat f, v_h \right>_{\Omega_{diff}^h} \\
			&\leq \lp \max_{i,\ell} \|\widetilde f - T^{k-1}_h f\|_{L^\infty(S^{i,\ell})} + \max_{i,\ell} \|\widehat f - T^{k-1}_h f\|_{L^\infty(S^{i,\ell})} \rp \int_{\Omega_{diff}^h} 1\cdot|v|dx \\
			&\leq C\delta_h^{k-1}|\Omega_{diff}^h|^{\frac12} |f|_{W^{k-1}_\infty(\Omega_h)}\|v_h\|_{0,\Omega_h} \\
			&\leq C\delta_h^{k-\frac12} |f|_{W^{k-1}_{\infty}(\Omega_h)} \|v_h\|_{1,\Omega_h},
		\end{aligned}
	$$
	after applying \cite[Lemma 4.3.8]{brenner2008mathematical}, H\"older's inequality, and seeing that $|\Omega^h_{diff}| \sim\mathcal O(\delta_h)$. Since we have assumed that $\delta_h \sim \mathcal O(h^2)$, we have that
	\begin{equation} \label{eqn: nonconformity inequality 1}
		\left< \widetilde f - \widehat f, v_h \right>_{\Omega_h} \leq Ch^{2k - 1}|f|_{W^{k-1}_\infty(\Omega_h)} \|v_h\|_{1,\Omega_h}.
	\end{equation}

	Next, using Lemma \ref{lemma: Taylor lemma 2}, we have that
	\begin{equation} \label{eqn: nonconformity inequality 2}
	\begin{aligned}
		\left<\widetilde p\circ\ex \res{\mathbf R^{k-1}_h(\nabla \widetilde u)}_\ex\cdot \nb,v_h \right>_{\Gamma_h}
		&\leq \overline p \lno\res{\mathbf R^{k-1}_h \nabla \widetilde u}_\ex\rno_{L^\infty(\Gamma_h)} \|v_h\|_{L^1(\Gamma_h)} \\
		&\leq C{\delta^k_h}|\nabla \widetilde u|_{W^{k}_\infty(\Omega)}\|v_h\|_{0,\Gamma_h} \\
		&\leq Ch^{2k}|u|_{W^{k+1}_\infty(\Omega)}\|v_h\|_{1,\Omega_h},
	\end{aligned}
	\end{equation}
	after applying the trace theorem and the assumption that $\delta_h \sim \mathcal O(h^2)$.

	Additionally, we have that
	\begin{equation} \label{eqn: nonconformity inequality 3}
	\begin{aligned}
		\left< \widetilde p\circ\ex \nabla e_h\cdot(\nb-\nb_h), v_h \right>_{\Gamma_h}
		&\leq \overline p\lp\max_{\xi \in \Gamma_h}\|\nb\circ\ex - \nb_h(\xib)\|_{\mathbb R^d}\rp \|\nabla e_h\|_{L^\infty(\Gamma_h)} \|v_h\|_{L^1(\Gamma_h)} \\
		&\leq Ch\|e_h\|_{W^{1}_\infty(\Omega_h)} \|v_h\|_{1,\Omega_h} \\
		&\leq Ch^{k+1} \lp |u|_{W^{k+1}_\infty(\Omega)} + |f|_{W^{k-1}_{\infty}(\Omega)} \rp \|v_h\|_{1,\Omega_h} ,
	\end{aligned}
	\end{equation}
	after applying Theorem \ref{theorem: max derivative optimality}.

	Then, we have that
	\begin{equation} \label{eqn: nonconformity inequality 4}
	\begin{aligned}
  		&\left< \widetilde p(\xib) \res{\Tb_h^{1,k-1}(\nabla e_h)}_\ex\cdot \nb, v_h\right>_{\Gamma_h} \\
  		&\quad\leq C\overline p\|v_h\|_{1,1\Omega_h}  \sum_{|\alphab|=1}^{k-1} \delta_h^{|\alphab|} \|D^\alphab \nabla e_h \|_{L^\infty(\Gamma_h)} \\
  		&\quad\leq C\overline p \|v_h\|_{1,1\Omega_h} \sum_{|\alphab|=1}^{k-1}\delta_h^{|\alphab|}  \lp\|D^\alphab \nabla(\widetilde u - u_I)\|_{L^\infty(\Omega_h)} + \|D^\alphab \nabla(u_I - u_h) \|_{L^\infty(\Omega_h)}\rp \\
  		&\quad\leq C\overline p\|v_h\|_{1,1\Omega_h}  \sum_{|\alphab|=1}^{k-1}\delta_h^{|\alphab|}  \lp h^{k-|\alphab|}|u|_{W^{k+1}_\infty(\Omega)} +Ch^{-|\alphab|} \lp \|\nabla(e_h)\|_{L^\infty(\Omega_h)} - \|\nabla(\widetilde u - u_I) \|_{L^\infty(\Omega_h)} \rp \rp \\
  		&\quad\leq C\overline p \|v_h\|_{1,1\Omega_h} \sum_{|\alphab|=1}^{k-1}\delta_h^{|\alphab|} h^{k-|\alphab|} \lp |u|_{W^{k+1}_\infty(\Omega)} + |f|_{W^{k-1}_\infty(\Omega)} \rp \\
  		&\quad\leq Ch^{k+1} \lp |u|_{W^{k+1}_\infty(\Omega)} + |f|_{W^{k-1}_\infty(\Omega)} \rp \|v_h\|_{1,\Omega_h} \\
	\end{aligned}
	\end{equation}
	after applying Theorem \ref{theorem: max derivative optimality} and the interpolation estimate \cite[Theorem 4.2.20]{brenner2008mathematical}.

	We finally move on to the last term, where we see that
	\begin{equation}\label{eqn: nonconformity inequality 5}
	\begin{aligned}
	\left< (\widetilde p(\xib) - \widetilde p\circ\ex)\res{\Tb^{k-1}_h(\nabla e_h)}_{\ex}\cdot\nb,v_h\right>_{\Gamma_h} 
	&\leq C\delta_h \lno \res{\Tb^{k-1}_h(\nabla e_h)}_{\ex} \rno_{L^\infty(\Gamma_h)} \|v_h \|_{L^1(\Gamma_h)} \\
	&\leq C\delta_h \|\nabla e_h \|_{L^\infty(\Omega_h)} \|v_h\|_{1,\Omega_h} \\
	&\leq Ch^{k+2}\lp |u|_{W^{k+1}_\infty(\Omega)} + |f|_{W^{k-1}_\infty(\Omega)} \rp \|v_h\|_{1,\Omega_h},
	\end{aligned}
	\end{equation}
	where we used Lemma \ref{lemma: Taylor lemma 1} and Theorem \ref{theorem: max derivative optimality}.

	Inserting \eqref{eqn: nonconformity inequality 1}, \eqref{eqn: nonconformity inequality 2}, \eqref{eqn: nonconformity inequality 3}, \eqref{eqn: nonconformity inequality 4}, and \eqref{eqn: nonconformity inequality 5} into \eqref{eqn: nonconformity inequality} gives us 
	$$
		N_h(\widetilde u - u_h, v_h) \leq Ch^{k+1}\lp |u|_{W^{k+1}_\infty(\Omega)} + |f|_{W^{k-1}_\infty(\Omega)} \rp \|v_h\|_{1,\Omega_h} \quad \forall v_h \in V^k_h
	$$
	This concludes this proof.
\end{proof}
Now that we have established that the nonconformity in the bilinear form is bounded above by $\mathcal O(h^{k+1})$, we are now ready to analyze the $L^2(\Omega)$ error of the numerical solution.

\subsection{The Dual Problem}
The dual problem we are interested in utilizing is to seek a $\psi \in H^1(\Omega_h)$ such that
$$
	N_h(v,\psi) = \left< \widetilde u - u_h, v \right>_{\Omega_h} \quad \forall v \in H^1(\Omega_h).
$$
The corresponding finite element approximation to the dual problem is to seek a $\psi_h \in V^k_h$ such that
$$
	N_h(v,\psi_h) = \left< \widetilde u-u_h,v\right>_{\Omega_h} \quad \forall v_h \in V^k_h.
$$
In general, $\phi \in H^1(\Omega_h)$ does not full $H^2(\Omega_h)$ regularity due to the presence of interior angles in $\Gamma_h$. From the literature, \cite{brenner2008mathematical, ciarlet2002finite} it is known that
\begin{equation} \label{eqn: dual bound}
	\|\psi\|_{1+s,\Omega_h} \leq C\|u-u_h\|_{0,\Omega_h}
\end{equation}
and
\begin{equation} \label{eqn: discrete dual bound}
	\|\phi\|_{1,\Omega_h} \leq C\|u-u_h\|_{-1,\Omega_h} \leq C\|u-u_h\|_{0,\Omega_h}.
\end{equation}

Furthermore, we have that
\begin{equation} \label{eqn: dual error}
	\|\psi - \psi_h\|_{1,\Omega_h} \leq Ch^s |\phi|_{1+s,\Omega_h},
\end{equation}
where $s\in\lp\frac12,1\right]$ depends on the magnitude of the largest interior angle of $\Gamma_h$. See \cite[Remark 7]{cheung2019optimally}.

\subsection{Analysis of $L^2(\Omega_h)$ Convergence}
With the previous results proven, we are now ready to derive optimal $L^2(\Omega_h)$ error estimates for the PE-FEM approximation of elliptic Neumann boundary value problems.
\begin{theorem}
	Assume that all the conditions in Theorem \ref{theorem: max derivative optimality} are satisfied. Then we have that
	$$
		\|\widetilde u - u_h \|_{0,\Omega_h} \leq Ch^{k+s}\lp |u|_{W^{k+1}_\infty(\Omega)} + |f|_{W^{k-1}_\infty\Omega} \rp,
	$$
	where $s\in\left(\frac12,1\right]$ depends on the largest interior angle of $\Gamma_h$. If $\Gamma_h$ is convex, then $s=1$.
\end{theorem}
\begin{proof}
From the definition of the dual problem, we have that
$$
	\begin{aligned}
		\lno \widetilde u - u_h \rno^2_{0,\Omega_h} &= N_h(\widetilde u - u_h, \phi) &\\
		& = N_h(\widetilde u - u_h, \phi - \phi_h) + N_h(\widetilde u - u_h, \phi_h) &\\
		&\leq C\|\widetilde u- u_h\|_{1,\Omega_h}\|\phi - \phi_h\|_{1,\Omega_h} + N_h(\widetilde u - u_h, \phi_h) &\textrm{(\cite[Theorem 2]{cheung2019optimally})} \\
		&\leq Ch^{k+s}\lp |u|_{k+1,\Omega_h} + |f|_{k-1,\Omega_h}\rp |\phi|_{2,\Omega_h}+ N_h(\widetilde u - u_h, \phi_h) &\textrm{\eqref{eqn: dual error} and \cite[Theorem 5]{cheung2019optimally}} \\
		&\leq Ch^{k+s}\lp |u|_{k+1,\Omega_h} + |f|_{k-1,\Omega_h}\rp |\phi|_{2,\Omega_h}& &\\
		&\quad + Ch^{k+1}\lp |u|_{W^{k+1}_\infty(\Omega)} + |f|_{W^{k-1}_\infty\Omega} \rp \|\phi_h\|_{1,\Omega_h} &\textrm{Lemma \ref{lemma: nonconformity error lemma}} \\
		&\leq Ch^{k+s}\lp |u|_{W^{k+1}_\infty(\Omega)} + |f|_{W^{k-1}_\infty\Omega} \rp \|\widetilde u - u_h \|_{0,\Omega_h} & \textrm{\eqref{eqn: dual bound} and \eqref{eqn: discrete dual bound}}.
	\end{aligned}
$$
This concludes this proof.
\end{proof}

\section{Conclusion}
In this manuscript, we have presented new results that indicate that the solution of the Neumann approximation presented in \cite{cheung2019optimally} is optimal in $W^1_\infty(\Omega)$ and $L^2(\Omega)$. Our analysis implies hat the solution to the continuous Neumann problem must be pointwise bounded a.e. in its $(k+1)$-th derivative in order for us to achieve additional accuracy in the $L^2(\Omega)$ norm, otherwise the $L^2(\Omega)$ error is bounded by the $\mathcal O(h^k)$ $H^1(\Omega)$ error. 

In our future work, we plan on utilizing these results in the derivation of optimal $L^2(\Omega)\times L^2(\Omega)$ error estimates for the interface coupling method presented in \cite{CHEUNG2020100094}. We will then aim to establish a generalized theory for the well-posedness and approximation properties of extension boundary methods. We believe that this approach can be generalized across many types of partial differential equations since the work presented in this manuscript and in \cite{cheung2019optimally} indicates that the averaged Taylor series approximation only generates an $\mathcal O(h)$ perturbation to the discretized variational operator that vanishes as $h\rightarrow 0$. 

\nocite{*}
\bibliographystyle{plain}
\bibliography{bibliography}

\begin{thebibliography}{1}

\bibitem{adams2003sobolev}
Robert~A Adams and John~JF Fournier.
\newblock {\em Sobolev spaces}.
\newblock Elsevier, 2003.

\bibitem{brenner2008mathematical}
Susanne~C Brenner, L~Ridgway Scott, and L~Ridgway Scott.
\newblock {\em The mathematical theory of finite element methods}, volume~3.
\newblock Springer, 2008.

\bibitem{CHEUNG2020100094}
James Cheung, Max Gunzburger, Pavel Bochev, and Mauro Perego.
\newblock An optimally convergent higher-order finite element coupling method
  for interface and domain decomposition problems.
\newblock {\em Results in Applied Mathematics}, 6:100094, 2020.

\bibitem{cheung2019optimally}
James Cheung, Mauro Perego, Pavel Bochev, and Max Gunzburger.
\newblock Optimally accurate higher-order finite element methods for polytopial
  approximations of domains with smooth boundaries.
\newblock {\em Mathematics of Computation}, 88(319):2187--2219, 2019.

\bibitem{ciarlet2002finite}
Philippe~G Ciarlet.
\newblock {\em The finite element method for elliptic problems}.
\newblock SIAM, 2002.

\end{thebibliography}
\end{document}